\documentclass[11pt]{elsarticle}
\usepackage{amssymb}
\usepackage{mathrsfs}
\usepackage{hyperref}
\usepackage{lmodern}

\newtheorem{thm}{Theorem}[section]
\newtheorem{cor}[thm]{Corollary}
\newtheorem{lem}[thm]{Lemma}

\newtheorem{prop}[thm]{Proposition}

\newcommand{\J}{\mathscr{J}}

\newcommand{\G}{\mathscr{G}}
\newcommand{\D}{\mathscr{D}}
\newcommand{\ran}{{\rm r.ann}_R}
\newcommand{\ann}{{\rm ann}_R}
\newcommand{\lan}{{\rm l.ann}_R}
\newproof{pf}{Proof}
\newcommand{\bx}{\widehat{x}}

\journal{Journal of Pure and Applied Algebra}

\begin{document}

\begin{frontmatter}

\title{{\bf On the semiprimitivity of free skew extensions of rings}}%

\author{Jeffrey Bergen\fnref{fn1}}%

\author[PG]{Piotr Grzeszczuk\corref{cor1}%
\fnref{fn2}}

\ead{p.grzeszczuk@pb.edu.pl}

\cortext[cor1]{Corresponding author}
\fntext[fn1]{Professor Jeffrey Bergen (DePaul University) passed away as this paper was close to be finished.}
\fntext[fn2]{The research of Piotr Grzeszczuk was supported by
the Bia{\l}ystok University of Technology grant WZ/WI/1/2019 and funded by the
resources for research by Ministry of Science and Higher Education of Poland.}

\address[PG]{Faculty of Computer Science, Bia{\l}ystok University of
Technology, Wiejska 45A, 15-351 Bia{\l}ystok, Poland}

\begin{abstract}
Let $X$ be a set of noncommuting variables, and $\G=\{\sigma_x\}_{x\in X}$, $\D=\{\delta_x\}_{x\in X}$ be sequences of automorphisms and skew derivations of a ring $R$. It is proved that if the ring $R$ is semiprime Goldie, then the free skew extension $R[X;\G,\D]$ is semiprimitive.
\end{abstract}

\end{frontmatter}

\section*{Introduction}

A well-known result of S. A. Amitsur \cite{A} states that if the ring $R$ has no nil ideals then the
polynomial ring $R[x]$ is semiprimitive. Subsequently, there has been a great deal of work examining the Jacobson radicals of more general ring extensions such as skew polynomial rings of automorphism type and of derivation type. 
For skew polynomial rings $R[x;\sigma]$ of automorphism
type, it was shown in \cite{BR} that even if $R$ is commutative and reduced, then the Jacobson radical  $\J(R[x;\sigma])$ can be nonzero.  Many authors, including C. R. Jordan
and D. A. Jordan \cite{JJ}, A. D. Bell \cite{B}, S. S. Bedi and J. Ram \cite{BR},\cite{R} and A. Moussavi \cite{M,MH} have extended Amitsur's result to skew polynomial rings of the form $R[x; \sigma, \delta$] and with certain additional conditions on $R$, where $\sigma$ is an automorphism (or monomorphism) of $R$, and $\delta$ is a $\sigma$-derivation of $R$. Very important and deep results on the Jacobson radical of skew polynomial rings of derivation type were obtained by A. Smoktunowicz (see \cite{S1}). For other results for such rings we refer to \cite{BMS}, \cite{BG1}, \cite{GSZ}, \cite{Mad} and \cite{SZ}.

Recall that for a given unital ring $R$ with a ring endomorphism  $\sigma$, a {\it $\sigma$-derivation} of $R$ is an additive map $\delta\colon  R\to R$  satisfying the $\sigma$-Leibniz rule
$$
\delta(xy)=\delta(x)y+\sigma(x)\delta(y),
$$ for $x,y\in R$. 
Then the skew polynomial ring $R[x;\sigma,\delta]$ can be described as the ring generated freely over $R$ by an element $x$ subject to the relation $xr=\sigma(r)x+\delta(r)$ for each $r\in R$. Throughout this paper we consider
skew polynomial rings over an arbitrary set of noncommuting variables. More precisely, let $X$ be a nonempty set and suppose that for any $x\in X$ there exists a pair $(\sigma_x,\delta_x)$, where $\sigma_x\colon R\to R$ is a ring endomorphism  and $\delta_x\colon R\to R$ is a  $\sigma_x$-derivation.
Put $\G=\{\sigma_x\}_ {x\in X}$ and $\D=\{\delta_x\}_ {x\in X}.$ Let us emphasize that we do not assume that the mappings $x\mapsto \sigma_x$ and $x\mapsto\delta_x$ are injective. Next we denote by $ \langle X\rangle$ the free monoid generated by the set of free generators $X$.  We will write $S=R[X;\G,\D]$ provided
\begin{enumerate}
\item $S$ contains $R$ as a subring;
\item $X\subseteq S$;
\item $S$ is a free left $R$-module with basis $\langle X\rangle$;
\item $xr=\sigma_x(r)x+\delta_x(r)$ for all $x\in X$ and $r\in R$.
\end{enumerate}
Such a ring $S$ we call a {\it free skew extension } of $R$. If $X=\{x\}$, $\mathscr{G}=\{\sigma\}$ and $\mathscr{D}=\{\delta\}$, then by Proposition 2.3 in \cite{GW} $S$ is the skew polynomial ring $R[x;\sigma,\delta]$. In Section 1 we prove that free skew extensions $S=R[X;\G,\D]$ always exist and they are unique up to isomorphism.
The main result of the paper, which we will prove in Section 2, is
\medskip

\noindent {\bf Theorem A.} {\em 
Let $R$ be a semiprime left Goldie ring. Then for any  set $X$ of noncommuting variables and sequences $\G=\{\sigma_x\}_{x\in X}$, $\D=\{\delta_x\}_{x\in X}$ of automorphisms and skew derivations of $R$ the free skew extension
 $R[X; \G,\D]$ is semiprimitive.}
\medskip

It is worth noting that A. D. Bell proved in \cite{B}  that if $\sigma$ is an automorphism and $\delta$ is a $\sigma$-derivation of a semiprime Goldie ring $R$, then $R[x;\sigma,\delta]$ is a semiprimitive Goldie ring. This result was then extended by A. Moussavi \cite{M}  and A. Moussavi, E. Hashemi \cite{MH} to the case when $\sigma$ is an injective endomorphism of $R$. Observe that if the set $X$ has at least two elements, then any free skew extension $S=R[X;\G,\D]$ has infinite Goldie dimension. Indeed, if $x,y\in X$, $x\neq y$, then the sum of left ideals $\sum\limits_{n=1}^\infty Sxy^n$ is direct, and so $S$ is not a Goldie ring.

Notice that if $\G=\{id_R\}_{x\in X}$, then $\D$ is a sequence of ordinary derivations of $R$. This special case was considered in \cite{TWC} under the name of the Ore extension of derivation type, and denoted by $R[X;\D]$. In Section 2 we apply Theorem A to free skew extensions of a certain class of rings, which includes rings with Krull dimension.
In particular we obtain

\medskip

\noindent {\bf Corollary B.} {\em Let $R$ be an algebra over a field $F$ of characteristic zero. If $R$ has left Krull dimension, then for any set $X$  of noncommuting variables and a sequence $\mathscr{D}=\{d_x\}_{x\in X}$ of derivations of $R$ the Jacobson radical  of the Ore extension $R[X;\mathscr{D}]$ is nilpotent. }

\section{Preliminary results and definitions}
In this Section  we prove the existence and uniqueness of free skew extensions. 
To do this, we will apply the arguments from the book by K. R Goodearl and R. B Warfield \cite{GW} used in the proof
of the existence and uniqueness of skew polynomial rings $R[x;\sigma,\delta]$.

\begin{prop} Given a unital ring $R$, a set $X$ of noncommuting variables and sequences $\G=\{\sigma_x\}_{x\in X}$, $\D=\{\delta_x\}_{x\in X}$ of ring endomorphisms and skew derivations of $R$, there exists a free skew extension 
$R[X;\G,\D]$.
\end{prop}

\begin{pf}
Let $R\langle X\rangle$ be a monoid ring of the free monoid $\langle X\rangle$ over $R$. Then any element of $R\langle X\rangle$ can be written uniquely in the form 
$$
f(X)=\sum\limits_{\Delta\in \langle X\rangle } r_\Delta \Delta \  \ ( {\rm  or\  shorter}\  \  f(X)=\sum r_\Delta \Delta),
$$
where $r_\Delta\in R$  and all but finitely many of these elements are zero.
Let $E={\rm End}_{\mathbb{Z}}(R\langle X\rangle)$ be the ring of all endomorphisms of the additive group $R\langle X\rangle$. Then $R$ can be treated as a subring of $E$ via left multiplications $r\mapsto r_{\ell}$; that is
$r_{\ell}(f(X))=rf(X),$ where $r\in R$, $f(X)\in R\langle X\rangle$.

For $x\in X$ define $\bx\in E$ according to the rule
$$
\bx(\sum r_\Delta \Delta)=\sum \left(\sigma_x(r_\Delta)x\Delta+\delta_x(r_\Delta)\Delta\right).
$$
Next, let $S$ be the subring of $E$ generated by $R\cup\{\bx\mid x\in X\}.$ Notice that for any $r\in R$ and any $f(X)=\sum r_\Delta \Delta$ in  $R\langle X\rangle$ 
$$
\begin{array}{rl}

\bx r(f(X))&=\bx(\sum rr_\Delta \Delta)=\sum(\sigma_x(rr_\Delta)x\Delta+\delta_x(rr_\Delta)\Delta)\\ \\

&= \sum(\sigma_x(r)\sigma_x(r_\Delta)x\Delta+\delta_x(r)r_\Delta\Delta+\sigma_x(r)\delta_x(r_\Delta)\Delta)\\ \\

&= \sigma_x(r)\sum(\sigma_x(r_\Delta)x\Delta+\delta_x(r_\Delta)\Delta) +\delta_x(r)\sum r_\Delta\Delta\\ \\

&= \sigma_x(r)\bx(f(X))+\delta_x(r)(f(X)=(\sigma_x(r)\bx+\delta_x(r))(f(X)).
\end{array}
$$
Thus 
\begin{equation}\label{relations}
\bx r=\sigma_x(r)\bx+\delta_x(r) \ \ {\rm for\    all\  } x\in X \ {\rm and\ } r\in R
\end{equation}
For $\Delta= x_1x_2\ldots x_n\in \langle X\rangle$ let $\widehat{\Delta}=\bx_1 \bx_2\ldots\bx_n$. The freeness of $\langle X\rangle$ implies that
$\langle\widehat{X}\rangle=\{\widehat{\Delta}\mid \Delta\in \langle X\rangle\}$ is a homomorphic image of the monoid
$\langle X\rangle$. The relations (\ref{relations}) imply that every element of $S$ can be written in the form $\sum r_{\widehat{\Delta}} \widehat{\Delta}$, so $S$  is generated by $\langle\widehat{X}\rangle$ as a left $R$-module. Observe that $\widehat{\Delta}(1)=\Delta$, so if
$\sum r_{\widehat{\Delta}} \widehat{\Delta}$ is a zero map, then 
$$
0=(\sum r_{\widehat{\Delta}} \widehat{\Delta})(1)=\sum r_{\widehat{\Delta}}\Delta.
$$
This forces that all coefficients $r_{\widehat{\Delta}}$  are equal to $0$, and hence $S$ is a free left $R$-module with basis 
$\langle\widehat{X}\rangle$. Notice also that the property 
$\widehat{\Delta}(1)=\Delta$, implies that the monoid $\langle\widehat{X}\rangle$ is free, and clearly isomorphic to $\langle X\rangle$. Thus the ring $S$ satisfies
all required conditions to be a free skew extension $R[X;\G,\D]$. \qed
\end{pf}

\medskip

Recall that the degree of a monomial $\Delta=x_1x_2\ldots x_n \in\langle X\rangle$ is defined to be $n$ and denoted by 
$\deg\Delta$.  
The uniqueness of  free skew extensions $R[X;\G,\D]$ will follow from the following universal mapping property

\begin{prop}\label{universal}
Let $S=R[X;\G,\D]$ be a free skew extension, where $\G=\{\sigma_x\}_{x\in X}$, $\D=\{\delta_x\}_{x\in X}$ are sequences of ring endomorphisms and skew derivations of $R$, respectively. Suppose that $T$ is a ring such that we have a ring homomorphism $\psi\colon R\to T$, and a mapping $\varphi\colon X\to T$ such that 
$$
\varphi(x)\psi(r)=\psi(\sigma_x(r))\varphi(x)+\psi(\delta_x(r))
$$ 
for all $x\in X$ and $r\in R$. Then there is a unique ring homomorphism
$\overline{\psi}\colon S\to T$ such that $\overline{\psi}_{| R}=\psi$ and $\overline{\psi}(x)=\varphi(x)$ for $x\in X$.
\end{prop}
\begin{pf}
Let $\overline{\varphi}\colon \langle X\rangle\to T$ be a unique extension of $\varphi$ to a homomorphism of monoids.
It is clear  that the only possibility for $\overline{\psi}$ is the map given by the rule
$$
\overline{\psi}(\sum r_\Delta\Delta)=\sum\psi(r_\Delta)\overline{\varphi}(\Delta)
$$

Notice that $\overline{\psi}$ is additive, and so it is enough to show that $\overline{\psi}(ab)=\overline{\psi}(a)\overline{\psi}(b)$ for $a,b\in S$.
Observe that if $b=\sum b_\Delta\Delta\in S$ and $x\in X$, then
$$
\begin{array}{rl}
\overline{\psi}(xb)&=\overline{\psi}\left(\sum\sigma_x(b_\Delta)x\Delta +\sum \delta_x(b_\Delta)\Delta\right)
=\sum \psi(\sigma_x(b_\Delta))\overline{\varphi}(x\Delta) +\sum \psi(\delta_x(b_\Delta))\overline{\varphi}(\Delta)\\ \\
&=\left( \sum \psi(\sigma_x(b_\Delta))\varphi(x) +\psi(\delta_x(b_\Delta))\right) \overline{\varphi}(\Delta)
= \sum \varphi(x)\psi(b_\Delta)\overline{\varphi}(\Delta)\\ \\ 
&= \varphi(x)\overline{\psi}(b)=\overline{\psi}(x)\overline{\psi}(b)
\end{array}
$$
Next it follows by induction on $\deg \Delta$  that $\overline{\psi}(\Delta b)=
\overline{\varphi}(\Delta)\overline{\psi}( b )$ for all $\Delta$ in $\langle X\rangle$. Furthermore if $r\in R$, then
$$
\overline{\psi}(rb)=\sum\psi(rb_\Delta)\overline{\varphi}(\Delta)=\sum \psi(r)\psi(b_\Delta)\overline{\varphi}(\Delta)=\psi(r)\sum \psi(b_\Delta)\overline{\varphi}(\Delta)
=\psi(r)\overline{\psi}(b).
$$
Consequently, if $a=\sum a_\Delta\Delta$, then
$$
\overline{\psi}(ab)=\sum\overline{\psi}(a_\Delta\Delta b)=\sum\psi(a_\Delta)\overline{\psi}(\Delta b)=
\sum \psi(a_\Delta)\overline{\varphi}(\Delta)\overline{\psi}(b)=\overline{\psi}(a)\overline{\psi}(b).
$$
This finishes the proof. \qed 
\end{pf}

As an immediate consequence of Proposition \ref{universal} we obtain the uniqueness of free skew extensions.

\begin{cor}
Let $R$ be a ring, $\G=\{\sigma_x\}_{x\in X}$, $\D=\{\delta_x\}_{x\in X}$ be sequences of ring endomorphisms and skew derivations of $R$, respectively. Suppose that $\varphi\colon X\to X^\prime$ is a bijection of sets, and put 
$\G^\prime=\{\sigma_{\varphi^{-1}(x\prime)}\}_{x^\prime\in X^\prime}$, $\D^\prime=\{\delta_{\varphi^{-1}(x^\prime)}\}_{x^\prime\in X^\prime}$
where $x^\prime=\varphi(x)$ for $x\in X$. If $S=R[X;\G,\D]$ and $S^\prime=R[X^\prime;\G^\prime,\D^\prime]$, then there exists a unique isomorphism $\psi\colon S\to S^\prime$ such that $\psi_{| R}$ is the identity map on $R$ and 
 $\psi(x)=x^\prime$ for $x\in X$.
\end{cor}

Let us fix a linear order $\prec$ on the set $X$. Then one can extend it lexicographically to the monoid $\langle X\rangle$; that is  $x_{i_1}\ldots x_{i_k}\prec x_{j_1}\ldots x_{i_l}$ if and only if
\begin{enumerate}
\item[$\bullet$] $k< l$ or
\item[$\bullet$] $k=l$ and there exists $s$, $1\leqslant s\leqslant k$ such that $i_1=j_1,\ldots, i_{s-1}=j_{s-1}$, $i_s\neq j_s$ and  $x_{i_s}\prec x_{j_s}$.
\end{enumerate}

If $0\neq f(X)=\sum r_{\Delta}\Delta\in R[X;\G,\D]$, then the set ${\rm supp} f(X)=\{\Delta\mid r_\Delta\neq 0\}$ is called the support of $f(X)$.
The largest monomial  $\Delta_0$ in the support of $f(X)$ (with respect to $\prec$) is called the leading term of $f(X)$. The coefficient $r_{\Delta_0}$ of the leading term is called the leading coefficient of $f(X)$. Thus any nonzero element $f=f(X)\in R[X;\G,\D]$ has a unique decomposition
\begin{equation}\label{decomp}
f= r\Delta_0+(f)_{\prec}
\end{equation}
where $\Delta_0$ is the leading term of $f$, $(f)_{\prec}$ is a finite sum of monomials $r_\Delta\Delta$ such that $\Delta\neq \Delta_0$ and $\Delta\prec\Delta_0$.

For a subset $A$ of $R$ we let $\lan(A)=\{r\in R\mid rA=0\}$, $\ran(A)=\{r\in R\mid Ar=0\}$ and $\ann(A)=\{r\in R\mid rA=Ar=0\}$ be the annihilators of $A$ in $R$. Recall that the ring $R$ is said to be left Goldie, if it contains no infinite direct sum of nonzero left ideals and satisfies  the ascending chain condition on left annihilators.

\section{Free skew extensions of semiprime Goldie rings}

Throughout this Section  we consider free skew extensions of automorphic type, that is we assume that $\G=\{\sigma_x\}_{x\in X}$ is a sequence of \underline{automorphisms} and $\D=\{\delta_x\}_{x\in X}$ is a sequence of respective
$\sigma_x$-derivations  of  a unital ring $R$.
Our goal is to prove that if the ring $R$ is semiprime Goldie, then any free skew extension $R[X;\G,\D]$ of automorphic type is semiprimitive. 

We start with the following general lemma.
\begin{lem}
Let $I$ be an ideal of a semiprime ring $R$ and let $c \in I$ such that $c$ is regular in $I$.
If $M = \ann(I)$, then
	\begin{enumerate}
		\item $\lan(c) = \ran(c) = M$ and
		\item $\eta(c)$ is regular in the factor ring $R/M$, where $\eta\colon R\to R/M$ is the canonical epimorphism $R$ onto $R/M$.
	\end{enumerate}
\end{lem}

\begin{pf}
We begin by observing that since $R$ is semiprime, the left and right annihilators of $I$ coincide
and have zero intersection with $I$.
Since $M = \ann(I)$, it is certainly the case that $Mc = cM = 0$.
Therefore to prove (1), we need to show that  $\lan(c)$ and $\ran(c)$ are contained in  $M$.
Suppose $r \in R$ such that $cr = 0$; therefore $crI = 0$.
Hence $rI$ consists of elements of $I$ which annihilate $c$ on the right.
However $c$ is regular in $I$, thus $rI = 0$.
As a result $r \in \ann(I) = M$.
An analogous argument shows that if $rc = 0$, then $r \in \ann(I) = M$.
Thus both $\lan(c)$ and $\ran(c)$ are contained in  $M$.
\vskip.1in

Next suppose $s \in R$ such that $cs \in M$; therefore $csI = 0$.
By the previous paragraph, we now have 
$sI \subseteq I \cap  \ran(c) \subseteq I \cap M = 0$.
Thus $s \in \ann(I) = M$.
This shows that the right annihilator of $\eta(c)$ in $R/M$ is zero and an analogous
argument works for the left annihilator. 
Thus $\eta(c)$ is regular in $R/M$.
\qed
\end{pf}

We continue with the following proposition.

\begin{prop}  Let $R$ be a semiprime left Goldie ring and let
$R[X;\G,\D]$ be a free skew extension of $R$.
If $A$ is a nonzero ideal of  $R[X;\G,\D]$  and $I $ is the additive subgroup of $R$
generated by the leading coefficients of elements of $A$ (with respect to a fixed linear order $\prec$ on $\langle X\rangle$),
then 
	\begin{enumerate}
		\item $I$ is a two-sided ideal of $R$,
		\item $\sigma_x(I) \subseteq I$, for every automorphism $\sigma_x$ corresponding to  $x \in X$
		\item if $M = \ann(I)$, then $\sigma_x(M) = M$ and $\delta_x(M) \subseteq M$, for every $x\in X$.
	\end{enumerate}
\end{prop}

\begin{pf}
For (1), it suffices to show that if $a \in I$ and $r \in R$, then $ra, ar \in I$.
Since $a \in I$,
$$
a = a_1 + \cdots + a_m,
$$
where each $a_i$ is the leading coefficient of some $\omega_i \in A$.
Therefore $ra_i$ is either $0$ or a leading coefficient of $r\omega_i \in A$.
Thus
$$
ra = ra_1 + \cdots + ra_m,
$$
hence $ra$ is the sum of leading coefficients of elements of $A$ and $ra \in I$.
\vskip.1in

According to (\ref{decomp}), each $\omega_i$ in the previous paragraph can be written as
$
\omega_i = a_i \Delta_i +  (\omega_i)_\prec,
$
where $\Delta_i$ is the leading term of $\omega_i$, and $(\omega_i)_\prec$ is a finite sum of monomials $r_\Delta\Delta$ such that $\Delta\neq\Delta_i$ and $\Delta\prec \Delta_i$.
If $s \in R$, then 
$$
\Delta_i s = \pi(s) \Delta_i + (\Delta_is)_\prec,
$$
for some automorphism $\pi$ of $R$. More precisely, if $\Delta=x_1x_2\ldots x_n$, then $\pi=\sigma_{x_1}\sigma_{x_2}\dots \sigma_{x_n}$. 
Therefore
$$
\omega_i \pi^{-1}(r) = a_ir \Delta_i +  (\omega_i \pi^{-1}(r))_\prec.
$$
Consequently, as above, $a_ir$ is either $0$ or the leading coefficient of $\omega_i \pi^{-1}(r)\in A$.
Thus
$$
ar = a_1r + \cdots + a_mr,
$$
hence $ar$ is the sum of leading coefficients of elements of $A$ and $ar \in I$.

\bigskip

For (2), let $a, a_i, \omega_i$ be as in the proof of (1).  Take $x\in X$.
Then $x\omega_i \in A$ and
$$
x\omega_i = \sigma_x(a_i) (x\Delta_i) +  (x\omega_i)_\prec
$$
Therefore $\sigma_x(a_i)$ is the leading coefficient of $x\omega_i \in A$.
Thus
$$
\sigma_x(a) = \sigma_x(a_1) + \cdots + \sigma_x(a_m),
$$
hence $\sigma_x(a)$ is the sum of leading coefficients of elements of $A$ and $\sigma_x(a) \in I$.

\bigskip

For part (3), observe that (2) gives us the following descending chain of two-sided ideals of $R$:
$$
I \supseteq \sigma_x(I)  \supseteq \sigma_x^2(I) \supseteq \sigma_x^3 (I) \supseteq  \cdots.
$$
Since $R$ is semiprime, if $i \geqslant 0$, the left and right annihilators of the ideal $\sigma_x^i(I)$ are the same.
In addition, since $R$ is left Goldie, it satisfies the descending chain condition on annihilators of ideals.
Thus there exists $n \in \mathbb{N}$ such that $\ann(\sigma_x^n(I)) =\ann(\sigma_x^{n+1}(I))$.
\vskip.1in

Whenever $\tau$ is an automorphism of $R$ and $C \subseteq R$, we have
$\tau(\ran(C)) = \ran(\tau(C))$.  
If $n \in \mathbb{N}$ is such that $\ann(\sigma_x^n(I)) = \ann(\sigma_x^{n+1}(I))$ and if we let $\sigma_x^n = \tau$,
then
$$\begin{array}{rl}
\ann(I) &= \sigma_x^{-n}(\sigma_x^n(\ann(I))) = \sigma_x^{-n}(\ann(\sigma_x^n(I))  = 
\sigma_x^{-n}(\ann(\sigma_x^{n+1}(I))\\ & =  \sigma_x^{-n}(\sigma_x^{n+1}(\ann(I)))  = \sigma_x(\ann(I)).
\end{array}
$$
The above equation shows that if $M = \ann(I)$, then $M =  \sigma_x(\ann(I))  = \sigma_x(M)$.
\vskip.1in

Take any $a,b \in I$ and $m\in M$. Then $\delta_x(a)b\in I$ and by (2) we obtain $\sigma_x(a) \delta_x(b)\in I$. Thus

$$
0 = \delta_x(abm) = \delta_x(a)bm + \sigma_x(a) \delta_x(b) m +  \sigma_x(a)   \sigma_x(b) \delta_x(m)  =\sigma_x(a)   \sigma_x(b) \delta_x(m).
$$
Therefore $(\sigma_x(I))^2 \delta_x(M) = 0$ and
$$
 (\sigma_x(I) \delta_x(M))^2 \subseteq   (\sigma_x(I))^2 \delta_x(M) = 0.
$$
Since $R$ is semiprime this tells us that the left ideal $\sigma_x(I) \delta_x(M)=0$,
hence
$$ 
\delta_x(M) \subseteq \ann(\sigma_x(I)) = \sigma_x(M) = M.
$$
Thus $\sigma_x(M) = M$ and $\delta_x(M) \subseteq M$.
\qed
\end{pf}

We will need the following lemma. 

\begin{lem}\label{disjoint} Let $X$ be a set of variables of cardinality $card(X)\geqslant 2$ and $\mathcal{F}=\{A_1,A_2,\ldots,A_m\}$ be a family of finite sets consisting of elements of the free monoid $\langle X\rangle$. Suppose that elements of $A_i$,
$i=1,2,\ldots, m$, have the same degree $n_i$. Then there exist an integer $t\geqslant \max\{n_i\mid 1\leqslant i\leqslant m\}$ and elements $\nu_1,\ldots,\nu_m\in \langle X\rangle$ such that $\deg\nu_i=t-n_i$ and the sets
$A_1\nu_1,A_2\nu_2,\ldots,A_m\nu_m$ are pairwise disjoint.
\end{lem}
\begin{pf}
The lemma is obvious when the set  $X$ is infinite. Suppose that $card(X)=d<\infty$. Without loss of generality we may assume that $n_1\geqslant n_2\geqslant\ldots\geqslant n_m$. It is clear that there are $d^s$ different elements of $\langle X\rangle$ of degree $s$. Take $s$ such that $d^s>m$ and choose different elements $\omega_1,\omega_2,\ldots,\omega_m\in \langle X\rangle$ of degree $s$. Fix an element $x\in X$ and put $\nu_i=x^{n_1-n_i}\omega_i$, $t=n_1+s$. Then  for any $a\in A_i$ $\deg a\nu_i=n_i+(n_1-n_i)+s=t$. Furthermore,  if $a\in A_i$, $b\in A_j$ and $i\neq j$, then 
$a\nu_i=(ax^{n_1-n_i})\omega_i\neq (bx^{n_1-n_j})\omega_j=b\nu_j$. This proves that the sets $A_1\nu_1,A_2\nu_2,\ldots,A_m\nu_m$ are pairwise disjoint.
\qed

\end{pf}

We can now prove our main result.

\begin{thm}\label{main}  Let $R$ be a semiprime left Goldie ring. Then for any  set $X$ of noncommuting variables and sequences $\G=\{\sigma_x\}_{x\in X}$, $\D=\{\delta_x\}_{x\in X}$ of automorphisms and skew derivations of $R$ the free skew extension
 $R[X; \G,\D]$ is semiprimitive.
\end{thm}

\begin{pf}
In light of the result of  A. D. Bell \cite{B} (mentioned in the Introduction) we may assume that the set $X$ contains at least two elements.
Let us fix a linear order $\prec$ on the free monoid $\langle X\rangle$.
By way of contradiction, let us suppose $\J(R[X; \G,\D]) \not= 0$
and we will apply Proposition 2 with $A = \J(R[X;\G,\D])$.
Next, we let $I$ be the ideal of leading coefficients in Proposition 2 and then let $M= \ann(I)$.
Recall that, for every $\sigma_x, \delta_x$ corresponding to some $x \in X$, we have
$\sigma_x(M) = M$ and $\delta_x(M) \subseteq M$.
\vskip.1in

Since $M= \ann(I)$, the sum $I + M$ is direct and is essential as both a left and right ideal of $R$.
Every essential one-sided ideal of a semiprime Goldie ring contains a regular element, therefore
there exist $c \in I$ and $d \in M$ such that $c + d$ is regular in $R$.
If $0 \not= s \in I$ then $sM = Ms = 0$, hence
$$
0 \not= s(c+d) = sc \  \  \   {\rm {and}}  \  \  \   0 \not= (c+d)s = cs.
$$
Thus $c$ is regular in $I$.  We can now apply Lemma 1 to $c \in I$ to conclude that
$\lan(c) = \ran(c) = M$ and $\eta(c)=c + M$ is regular in the factor ring $R/M$.
\vskip.1in

Since $0 \not=  c \in I$, there exist
$a_1, \ldots , a_m$ such that 
$c = a_1 + \cdots + a_m$, where each $a_i$ is a leading coefficient of some $\omega_i \in \J(R[X;\G,\D])$.
Using Lemma \ref{disjoint} we can multiply each $\omega_i$ by an appropriate monomial $\nu_i$ such that
all the $\omega_i\nu_i$ have the same degree, say $t$, but none of the monomials of highest degree in the support of $\omega_i\nu_i$ appear in the support of
$w_j\nu_j$, for $i \not= j$. According to (\ref{decomp}) we have 
$$\omega_i\nu_i=a_i\Delta_i+(\omega_i\nu_i)_\prec, \  \  
\deg\Delta_i=t \  \  {\rm and}\  \  \Delta_i\neq\Delta_j \  \  {\rm for}\  \ i\neq j.
$$
Therefore, if
$$
a(X) = \omega_1\nu_1 + \cdots + \omega_m\nu_m,
$$
we have $a(X) \in \J(R[X;\G,\D])$ and 
$$
a(X) = a_1\Delta_1 + \cdots + a_m\Delta_m +   \sum_i (\omega_i\nu_i)_\prec.
$$
Since none of the monomials $\Delta_i$ appear in the support of $\sum\limits_i (\omega_i\nu_i)_\prec$, we have
an irreducible decomposition
$$
a(X) = a_1\Delta_1 + \cdots + a_m\Delta_m +  f(X),
$$
where $\deg f(X)\leqslant t$.
In addition, since $a(X)x \in \J(R[X;\G,\D])$ for all $x \in X$, we may assume that the $\nu_i$ were chosen to make the degree of $a(X)$ equal to some  $t \geqslant 1$ and the constant term of $a(X)$ equal  to $0$.
\vskip.1in

Every element of $\J(R[X;\G,\D])$ is quasi-invertible, therefore there exists an element
$b(X) \in \J(R[X;\G,\D])$ such that
\begin{equation}\label{Eq1}
a(X) + b(X) = a(X) b(X) = b(X) a(X).
\end{equation}
\vskip.1in

For every $\sigma_x, \delta_x$ corresponding to some $x \in X$, we have $\sigma_x(M) = M$ and $\delta_x(M) \subseteq M$, 
therefore the actions of $\sigma_x$ and $\delta_x$  induce actions on the factor $R/M$.  
In addition, we can examine the free skew extension $M[X;\G,\D]$,
and observe that $M[X;\G,\D]$ is a two-sided ideal of $R[X;\G,\D]$.
We can now identify the factor ring 
$R[X;\G,\D]/M[X;\G,\D]$ with the free skew extension $(R/M)[X;\overline{\G},\overline{\D}]$.
If we let $\overline{a(X)}$ and $\overline{b(X)}$ be the images of $a(X)$ and $b(X)$ in $(R/M)[X;\overline{\G},\overline{\D}]$,
then equation (\ref{Eq1}) becomes
\begin{equation}\label{Eq2}
\overline{a(X)} + \overline{b(X)} = \overline{a(X)}  \cdot \overline{b(X)} = \overline{b(X)} \cdot  \overline{a(X)}.
\end{equation}
\vskip.1in

Recall that $c \notin M$, therefore at least one of $a_1, \ldots, a_m$ is not in $M$, 
hence $ \overline{a(X)}$ also has degree $t \geqslant 1$.
Thus equation (\ref{Eq2}) immediately implies that $\overline{b(X)}$ is not equal to zero in $(R/M)[X;\overline{\G},\overline{\D}]$. 
Now suppose $\overline{b(X)}$ has degree at least one.  Then there exists $b^\prime \in R$ 
and a monomial $\Delta$ of degree at least one 
such that $b^\prime \notin M$ and $\overline{b^\prime}\Delta$ is the leading term
of $\overline{b(X)}$.
Therefore 
$$
\overline{b(X)}  \cdot \overline{a(X)} = \overline{b^\prime \pi(a_1)} \Delta \Delta_1 + \cdots + \overline{b^\prime \pi(a_m)} \Delta \Delta_m +  \overline{g(X)},
$$
where $\pi$ is the automorphism of $R$ equal to the product of the automorphisms corresponding to the $x \in X$ appearing in $\Delta$, and non of the monomials $\Delta\Delta_i$ is contained in the support of $\overline{g(X)}\in (R/M)[X;\overline{\G},\overline{\D}]$.
Observe that, for $1 \leqslant i \leqslant m$, $\overline{b^\prime \pi(a_i)}$ must be the coefficient of $\Delta \Delta_i$ in $\overline{b(X)}  \cdot \overline{a(X)}$
since no other product of monomials from $\overline{a(X)} $ and $\overline{b(X)}$  could result in $\Delta \Delta_i$.
If any of the $\overline{b^\prime \pi(a_i)}$ is nonzero in $R/M$, then the degree of 
$\overline{b(X)}  \cdot \overline{a(X)}$ exceeds the degree of $\overline{a(X)} + \overline{b(X)}$, contradicting equation (\ref{Eq2}).
\vskip.1in

As a result, $b^\prime \pi(a_i) \in M$, for $1 \leqslant i \leqslant m$.
This implies that
$$
b^\prime\pi(c) = b^\prime \pi(a_1)  + \cdots + b^\prime \pi(a_m) \in M.
$$
Since $\pi(M) = M$, we have
$$
\pi^{-1}(b^\prime)c = \pi^{-1}(b^\prime\pi(c)) \in M.
$$
However, $\eta(c)=c + M$ is regular in $R/M$, hence $\pi^{-1}(b^\prime) \in M$.
This immediately implies that $b^\prime \in M$, contradiction that $b^\prime \not\in M$.
\vskip.1in

Having shown that $\overline{b(X)}$ is nonzero and had degree less than $1$ in $(R/M)[X;\overline{\G},\overline{\D}]$,
it follows that 
$$
cb(X) = cb_0 \not= 0,
$$
for some $b_0 \in R$.

Multiplying equation (\ref{Eq1}) on the left by $c$ now gives us
$$
ca(X) + cb(X) = ca(X) b(X) =c b(X) a(X).
$$
Since $cb(X) = c b_0$, examining the equation $ca(X) + cb(X) = cb(X)c(x)$ 
gives us
$$
ca(X) + cb_0 = cb_0 a(X).
$$

Since the constant term of $a(X)$ is $0$, if we compare the constant terms of both sides of the previous equation, we obtain
$cb_0 = 0$,  contradicting that $cb_0 \not= 0$. 
Consequently, the ring $R[X;\G,\D]$ is semiprimitive.
\qed
\end{pf}
\medskip

It is well known that any ring $R$ with left Krull dimension  has a nilpotent prime radical $\mathscr{P}(R)$  (cf. \cite{MR}, Corollary 6.3.8), and the factor ring $R/\mathscr{P}(R)$ is left Goldie (cf. \cite{MR}, Proposition 6.3.5). These properties of rings with Krull dimension motivate the following general observation

\begin{cor}\label{cor} Let $\G=\{\sigma_x\}_ {x\in X}$ and $\D=\{\delta_x\}_ {x\in X}$ be families of automorphisms  and skew derivations of the ring $R$ such that
\begin{enumerate}
\item the prime radical ${\mathscr{P}}(R)$ is nilpotent,
\item the factor ring $R/{\mathscr{P}}(R)$ is left Goldie,
\item ${\mathscr{P}}(R)$ is stable under $\delta_x$, for all $x\in X$.
\end{enumerate}
Then the Jacobson radical $\J(R[X;\G,\D])$ is nilpotent.
\end{cor}
\begin{pf}
Since ${\mathscr{P}}(R)$ is stable under all $\sigma_x$ and $\delta_x$, ${\mathscr{P}}(R)[X;\G,\D]$ is a nilpotent ideal of 
the skew extension $R[X;\G,\D]$. In particular, ${\mathscr{P}}(R)[X;\G,\D]\subseteq \J(R[X;\G,\D])$. Observe that 
$$
R[X;\G,\D]/{\mathscr{P}}(R)[X;\G,\D] \simeq
(R/{\mathscr{P}}(R))[X;\overline{\G},\overline{\D}],
$$ 
where $\bar{\sigma}_x$, $\bar{\delta}_x$ are induced automorphisms and skew derivations of $R/{\mathscr{P}}(R)$. The ring $R/{\mathscr{P}}(R)$ is semiprime Goldie, so by Theorem \ref{main} the free skew extension $(R/{\mathscr{P}}(R))[X;\overline{\G},\overline{\D}]$ is semiprimitive. It means that
$\J(R[X;\G,\D])={\mathscr{P}}(R)[X;\G,\D]$. 
\qed
\end{pf}

The problem when the nil and prime radicals of a ring  are stable under skew derivations is examined in \cite{BG2}. In particular, Lemma 3 of \cite{BG2} states that if $\delta$ is a $\sigma$-derivation of a ring $R$, then the nil radical 
$\mathscr{N}(R)$ is $\delta$-stable provided the element $\delta(a)$  is nilpotent for any $a\in \mathscr{N}(R)$. 
Notice that the condition (2) of Corollary \ref{cor} implies immediately the equality of radicals  $\mathscr{N}(R)=\mathscr{P}(R)$. Indeed, it is well known that  nil ideals of rings with ascending chain condition on left annihilators contain nonzero nilpotent ideals (see Lemma 2.3.2 of \cite{MR}). As a consequence, the condition (3) of Corollary \ref{cor} can be raplaced by:
\medskip
\begin{enumerate}
\item[($3^\prime$)] $\delta_x(a)$ is nilpotent for all $x\in X$ and $a\in \mathscr{P}(R)$. 
\end{enumerate}
\medskip

For the remainder of this paper, we will examine algebras over a field $F$ with  $q$-skew derivations. Recall that a $\sigma$-derivation $\delta$ of $R$ is said to be a $q$-skew $\sigma$-derivation if there exists a nonzero element $q\in F$ such that $\delta\sigma=q\sigma\delta$. From  Lemma 4 of \cite{BG2}  it follows that  if $I$ is a $\sigma$-stable ideal of $R$, then for any $a_1,a_2,\ldots,a_n\in I$ 
\begin{equation}\label{prod}
\delta^n(a_1a_2\ldots a_n)= (n!)_q \sigma^{n-1}\delta(a_1)\sigma^{n-2}\delta(a_2)\ldots \sigma\delta(a_{n-1})\delta(a_n) +w,
\end{equation}
where $w\in I$ and $(n!)_q=\prod\limits_{i=1}^n(1+q+\ldots+q^{i-1})$. Thus if $I^n=0$, then the identity (\ref{prod})
gives that for any $a\in I$ $(n!)_q(\delta(a))^n\in I$, and hence $((n!)_q)^n(\delta(a))^{n^2}=0$. As a consequence we obtain that $\delta$ satisfies ($3^\prime$) provided $(n!)_q\neq 0$ in $F$. Notice that $(n!)_q\neq 0$ in $F$ means that either $q$  is not a root o unity of degree $d\leqslant n$ or $n< {\rm char}\,F$, when $q=1$. Consequently, we obtain the following observation.

\begin{cor} Let $\G=\{\sigma_x\}_ {x\in X}$ and $\D=\{\delta_x\}_ {x\in X}$ be families of automorphisms  and skew derivations of an $F$-algebra $R$ such that
\begin{enumerate}
\item the prime radical ${\mathscr{P}}(R)$ is nilpotent with index of nilpotency equal to $n$,
\item the factor ring $R/{\mathscr{P}}(R)$ is left Goldie,
\item for any $x\in X$  $\delta_x\sigma_x=q_x\sigma_x\delta_x$, where $q_x\in F^*$ and $(n!)_{q_x}\neq 0$ in $F$.
\end{enumerate}
Then the Jacobson radical $\J(R[X;\G,\D])$ is nilpotent.
\end{cor}

Since the nil and prime radicals of algebras over fields of characteristic zero are stable under ordinary derivations (see also Proposition 2.6.28 in \cite{RL}), we obtain
 
\begin{cor} Let $R$ be an algebra over a field $F$ of characteristic zero. If $R$ has left Krull dimension, then for any set $X$  of noncommuting variables and a sequence $\mathscr{D}=\{d_x\}_{x\in X}$ of derivations of $R$ the Jacobson radical  of the Ore extension $R[X;\mathscr{D}]$ is nilpotent.
\end{cor}

\section*{Acknowledgements} The author would like to thank the anonymous referee for numerous helpful comments and suggestions, which have significantly improved this paper.


\end{document}